\documentclass[11pt]{article}
\usepackage{amsmath, amssymb, amscd, epsfig}

\newtheorem{thm}{Theorem}
\newtheorem{prop}[thm]{Proposition}
\newtheorem{lem}[thm]{Lemma}

\newtheorem{cor}[thm]{Corollary}

\renewcommand{\epsilon}{\varepsilon}
\renewcommand{\phi}{\varphi}
\renewcommand{\deg}{\widetilde{\operatorname{deg}}}

\newcommand{\BB}{\mathbb}

\newcommand{\pf}{\noindent {\it Proof. }}
\newcommand{\qed}{\nopagebreak $\qquad$ $\square$ \vskip5pt}
\newcommand{\separate}{\vskip5pt}
\newcommand{\supp}{\operatorname{supp}}
\newcommand{\im}{\operatorname{Im}}
\newcommand{\re}{\operatorname{Re}}

\newcommand{\B}{\overline}
\newcommand{\HC}{\BB H_{\BB C}}
\newcommand{\HR}{\BB H_{\BB R}}

\newcommand{\Pl}{\operatorname{Pl}}
\newcommand{\FC}{{\cal F}'_{C^*}}

\begin{document}

\title{\bf Quaternionic Analysis and the Schr\"odinger Model
for the Minimal Representation of $O(3,3)$}
\author{Igor Frenkel and Matvei Libine}
\maketitle

\begin{abstract}
In the series of papers \cite{FL,FL2} we approach quaternionic analysis from
the point of view of representation theory of the conformal group
$SL(2,\HC) \simeq SL(4,\BB C)$ and its real forms.
This approach has proven very fruitful and pushed further
the parallel with complex analysis and develop a rich theory.
In \cite{FL2} we study the counterparts of Cauchy-Fueter and Poisson formulas
on the spaces of split quaternions $\HR$ and Minkowski space $\BB M$ and
show that they solve the problem of separation of the discrete and continuous
series on $SL(2,\BB R)$ and the imaginary Lobachevski space
$SL(2,\BB C)/SL(2,\BB R)$.
In particular, we introduce an operator $\Pl_R$, compute its effect on
the discrete and continuous series components of the space of functions
${\cal H}(\HR^+)$ and obtain a surprising formula for the Plancherel
measure of $SL(2,\BB R)$.
The proof is based on a transition to the Minkowski space
$\BB M$ and some pretty lengthy computations.
In this paper we introduce an operator $\frac{d}{dR} \Pl_R$ on
${\cal H}(\HR^+)$ and show that its effect on the discrete and continuous
series components can be easily computed using the Schr\"odinger model for
the minimal representation of $O(p,q)$ (with $p=q=3$) and the results of
Kobayashi-Mano from \cite{KobM}, particularly their computation of the integral
expression for the operator ${\cal F}_C$.
This provides an independent verification of the coefficients involved
in the formula for $\Pl_R$.
This paper once again demonstrates a close connection between quaternionic
analysis and representation theory of various $O(p,q)$'s.
\end{abstract}

\section{Introduction}

In the series of papers \cite{FL,FL2} we approach quaternionic analysis from
the point of view of representation theory of the conformal group
$SL(2,\HC) \simeq SL(4,\BB C)$ and its real forms such as
$SL(2, \BB H) \approx SO(5,1)$, $SU(2,2) \approx SO(4,2)$,
$SL(4,\BB R) \approx SO(3,3)$.
This approach has proven very fruitful and pushed further
the parallel with complex analysis and develop a rich theory.
In \cite{FL2} we study the counterparts of Cauchy-Fueter and Poisson formulas
on the spaces of split quaternions $\HR$ and Minkowski space $\BB M$ and
show that they solve the problem of separation of the discrete and continuous
series on $SL(2,\BB R)$ and the imaginary Lobachevski space
$SL(2,\BB C)/SL(2,\BB R)$.
The continuous series component on $SL(2,\BB R)$ gives rise to the minimal
representation of the conformal group $SL(4,\BB R)$.
We also obtain a surprising formula for the Plancherel measure of $SL(2,\BB R)$
in terms of the Poisson integral on the split quaternions $\HR$ (Theorem 94).
For convenience we restate it here.

\begin{thm}  \label{cont_ser_proj}
The operator defined on the space of functions ${\cal H}(\HR^+)$
$$
\phi(X) \mapsto (\Pl_R \phi)(W) =
\lim_{\epsilon \to 0^+} \frac1{\pi^2} \int_{X \in H_R}
\biggl( \frac {1}{N(X-W) +i\epsilon} - \frac {1}{N(X-W) -i\epsilon} \biggr)
\phi(X) \,\frac{dS}{\|X\|}
$$
annihilates the discrete series component of ${\cal H}(\HR^+)$.
If $\re l =-\frac12$ and $W \in \HR^+$,
$$
(\Pl_R t^l_{n\,\underline{m}})(W) =
\bigl(1 + R^{4 \im l} \cdot N(W)^{-2l-1} \bigr) \cdot t^l_{n\,\underline{m}}(W) \cdot
\begin{cases}
\frac{\coth(\pi\im l)}{\im l} & \text{if $m,n \in \BB Z$;}  \\
\frac{\tanh(\pi\im l)}{\im l} & \text{if $m,n \in \BB Z +\frac12$.}
\end{cases}
$$
\end{thm}

Of course, strictly speaking, the restrictions to $SU(1,1)$ of the matrix
coefficient functions of the continuous series representations
$t^l_{n\,\underline{m}}(X)$, $l=-1/2+i\lambda$ with $\lambda \in \BB R$,
do not belong to $L^2(SU(1,1))$, so we need to work with wave packets,
as explained in Remark 93 in \cite{FL2}.

It is relatively easy to see that the operator $\Pl_R$
annihilates the discrete series and that $\Pl_R$ is
$SU(1,1) \times SU(1,1)$-equivariant, but it is much harder to determine
its effect on the continuous series.
The proof given in \cite{FL2} is based on a transition to the Minkowski space
$\BB M$ and some pretty lengthy computations.
It is highly desirable to find a more simple or at the very least a less
computational proof.

In this paper we give an elementary proof of a slightly weaker result
(Theorem \ref{cont_ser_proj-2}). This proof uses the setting of the
Schr\"odinger model for the minimal representation\footnote
{Strictly speaking, $SO(3,3)$ does not have a minimal representation,
so by ``minimal representation of $SO(3,3)$'' we mean the representation
$(\varpi^{p,q}, V^{p,q})$ of $O(p,q)$ in the notations of \cite{KobO},
\cite{KobM} for $p=q=3$.
When $p+q$ is an even number greater than or equal $8$, one gets a genuine
minimal representation.} of $O(p,q)$
(with $p=q=3$) and the results of Kobayashi-Mano from \cite{KobM},
particularly their computation of the integral expression for the operator
${\cal F}_C$. This provides an independent verification of the coefficients
involved in Theorem \ref{cont_ser_proj}.
It is interesting to note that our argument can be reversed and, assuming
Theorem \ref{cont_ser_proj-2}, one can obtain the integral expression
for the operator ${\cal F}_C$ when $p=q=3$.
The Minkowski space $\BB M$ is also implicitly present in this article.
We are essentially computing the ratio of two integral kernels from \cite{KobM}
-- one for $\HR$ and the other for $\BB M$.

In addition to proving Theorem \ref{cont_ser_proj-2}, we find all $K$-types of
the minimal representation of $SO(3,3)$
(in the Schr\"odinger model realization).
Some -- but not all -- $K$-types of the minimal representations $V^{p,q}$
of various $O(p,q)$'s were found in \cite{KobM}.
Kobayashi and Mano called the $K$-types they found the ``skeleton'' of the
minimal representation in the sense that they contain at least one non-zero
vector from each $K$-irreducible component of $V^{p,q}$.
Moreover, as was pointed out by the referee, more $K$-types were found by
Hilgert, Kobayashi, Mano and M\"ollers in \cite{HKMM}.
More precisely, they fix a certain subgroup $K'$ of $K= O(p) \times O(q)$
isomorphic to $O(p-1) \times O(q-1)$;
then Corollary 8.2 in \cite{HKMM} describes all $K$-types of $V^{p,q}$
that are fixed by $K'$ and asserts that each $K$-irreducible component
of $V^{p,q}$ contains a unique (up to scaling) non-zero $K'$-fixed vector.

Our paper once again demonstrates a close connection between quaternionic
analysis and representation theory of various $O(p,q)$'s.

The first author was supported by the NSF grants DMS-0457444 and DMS-1001633;
the second author was supported by the NSF grant DMS-0904612.

\section{Statement of the Main Result}

Let us recall the $\delta$-functions on the cone $C=\{X \in \HR ;\: N(X)=0 \}$
and the hyperboloids $H_R =\{X \in \HR ;\: N(X)=R^2 \}$, denoted by
$\delta(C)$ and $\delta(H_R)$ respectively,
$$
\delta(C) = \frac1{2\pi i} \biggl( \frac1{N(X)-i0} - \frac1{N(X)+i0} \biggr):
\quad \psi \mapsto \frac 12 \int_{X \in C} \psi(X) \,\frac{dS}{\|X\|},
$$
$$
\delta(H_R) =
\frac1{2\pi i} \biggl( \frac1{N(X)-R^2-i0} - \frac1{N(X)-R^2+i0} \biggr):
\quad \psi \mapsto \frac 12 \int_{X \in H_R} \psi(X) \,\frac{dS}{\|X\|}.
$$
Using these distributions we can formally rewrite the operator $\Pl_R$
introduced in \cite{FL2} as
$$
\Pl_R: \quad \phi \mapsto
\frac4{\pi i} \bigl( \phi \cdot \delta(H_R) \bigr) \ast \delta(C),
$$
where $\ast$ denotes the convolution and either $\delta(C)$ or $\delta(H_R)$
should be replaced with
\begin{align*}
\delta(C)(\psi) &= \lim_{\epsilon \to 0^+} \frac1{2\pi i} \int_{X \in \HR}
\biggl( \frac1{N(X)-i\epsilon} - \frac1{N(X)+i\epsilon} \biggr)\psi(X) \,dV
\qquad \text{or}\\
\delta(H_R)(\psi) &= \lim_{\epsilon \to 0^+} \frac1{2\pi i} \int_{X \in \HR} \biggl(
\frac1{N(X)-R^2-i\epsilon} - \frac1{N(X)-R^2+i\epsilon} \biggr)\psi(X) \,dV,
\end{align*}
with $dV= dx^0 \wedge dx^1 \wedge dx^2 \wedge dx^3$.
Thus we can rewrite  $\Pl_R$ as
$$
(\Pl_R \phi)(W) =
\lim_{\epsilon \to 0^+} \frac1{\pi^2} \int_{X \in C} \biggl(
\frac1{N(W-X)-R^2 +i\epsilon} - \frac1{N(W-X)-R^2 -i\epsilon} \biggr)
\phi(W-X) \,\frac{dS}{\|X\|}.
$$
Differentiating with respect to $R^2$ we obtain the following result:
\begin{thm}  \label{cont_ser_proj-2}
The operator
\begin{multline*}
\phi(X) \mapsto (\Pl'_R \phi)(W) \\
= \lim_{\epsilon \to 0^+} \frac1{2\pi^2} \int_{X \in C} \biggl(
\frac1{\bigl(N(W-X)-R^2 +i\epsilon\bigr)^2} -
\frac1{\bigl(N(W-X)-R^2 -i\epsilon\bigr)^2} \biggr) \phi(W-X) \,\frac{dS}{\|X\|}
\end{multline*}
annihilates the discrete series component of ${\cal H}(\HR^+)$.
If $\re l =-\frac12$ and $W \in \HR^+$,
$$
(\Pl'_R t^l_{n\,\underline{m}})(W) =
R^{4l} \cdot N(W)^{-2l-1} \cdot t^l_{n\,\underline{m}}(W) \cdot
\begin{cases}
\coth(\pi\im l) & \text{if $m,n \in \BB Z$;}  \\
\tanh(\pi\im l) & \text{if $m,n \in \BB Z +\frac12$.}
\end{cases}
$$
\end{thm}

The point of introducing the operator $\Pl'_R$ is that its Fourier transform
is particularly easy to compute.

\section{Schr\"odinger Model for the Minimal Representation of $O(3,3)$}

Let $w_0 = \begin{pmatrix} I_p & 0 \\ 0 & -I_q \end{pmatrix} \in O(p,q)$,
where $I_p$ and $I_q$ denote the $p \times p$ and $q \times q$ identity
matrices.
One of the main results of \cite{KobM} is finding an explicit expression
for the action of $w_0$ on the minimal representation of $O(p,q)$ realized
in the Schr\"odinger model $(\pi, L^2(C^*))$.
This explicit expression is then used to prove unitarity of the minimal
representation. In this paper we are interested in the case $p=q=3$.
Let $w'_0 \in O(3,3)$ be the diagonal matrix with diagonal entries
$(-1,1,1,1,1,1)$. The action of $O(3,3)$ on $\BB R^{2,2}$ is described in
Subsection 2.8 of Part III of \cite{KobO}, and it is easy to see that the
element $w'_0$ acts on $\BB R^{2,2}$ by
$$
w'_0 : \quad (u',u'') \mapsto \frac{4}{|u'|^2-|u''|^2} (u',u'').
$$
This action lifts to the solutions of the ultrahyperbolic wave equation:
$$
w'_0 : \quad \phi(u',u'') \mapsto  \frac4{|u'|^2-|u''|^2} \cdot
\phi \biggl( \frac{4}{|u'|^2-|u''|^2} (u',u'') \biggr).
$$
We identify the split quaternions $\HR$ with $\BB R^{2,2}$ as follows:
$$
\HR = \biggl\{ X = x_1 e_0 + x_3 \tilde e_1 + x_4 \tilde e_2 + x_2 e_3 =
\begin{pmatrix} x_1-ix_2 & x_3+ix_4 \\ x_3-ix_4 & x_1+ix_2 \end{pmatrix}
;\: x_1, x_2, x_3, x_4 \in \BB R \biggr\},
$$
$$
\HR \ni X \longleftrightarrow (u', u'') \in \BB R^{2,2},
\qquad u'=(x_1,x_2), \quad u''=(x_3,x_4),
$$
this way
$$
\square_{2,2} = \frac {\partial^2}{\partial x_1^2} +
\frac {\partial^2}{\partial x_2^2} -
\frac {\partial^2}{\partial x_3^2} -
\frac {\partial^2}{\partial x_4^2}.
$$
Then $w'_0$ corresponds to the operator on $\HR$
and the solutions of $\square_{2,2} \phi=0$
$$
w'_0 : \quad X \mapsto \frac{4X}{N(X)}, \qquad
\phi(X) \mapsto \frac4{N(X)} \cdot \phi \biggl( \frac{4X}{N(X)} \biggr).
$$
In particular, $w'_0$ acts on the matrix coefficients $t^l_{n\,\underline{m}}$'s
of $SU(1,1)$ (see Subsection 2.5 in \cite{FL2}) by
\begin{equation}  \label{w_0-action}
w'_0 : \quad t^l_{n\,\underline{m}}(X) \mapsto
2^{4l+2} \cdot N(X)^{-2l-1} \cdot t^l_{n\,\underline{m}}(X).
\end{equation}

As in \cite{KobM}, we use the following normalization of the Fourier transform
on $\HR$:
\begin{equation}  \label{FT-normalization}
\hat \psi(\xi) =
\frac1{(2\pi)^2} \int_{X \in \HR} \psi(X) \cdot e^{i \xi \cdot X} \,dV,
\end{equation}
where $\xi=(\xi_1,\xi_2,\xi_3,\xi_4) \in \HR^*$ -- the dual of $\HR$ and
$\xi \cdot X = \xi_1x_1+\xi_2x_2+\xi_3x_3+\xi_4x_4$.
For $\xi, \xi' \in \HR^*$, write
$$
\langle \xi,\xi' \rangle = \xi_1\xi'_1+\xi_2\xi'_2-\xi_3\xi'_3-\xi_4\xi'_4,
\qquad \xi=(\xi_1,\xi_2,\xi_3,\xi_4), \quad \xi'=(\xi'_1,\xi'_2,\xi'_3,\xi'_4).
$$
Let $C^*$ denote the cone dual to $C$:
$$
C^* = \bigl\{ \xi \in \HR^* ;\:
\langle \xi,\xi \rangle = (\xi_1)^2+(\xi_2)^2-(\xi_3)^2-(\xi_4)^2=0 \bigr\},
$$
and let ${\cal F}_{C^*}$, $\FC$ be integral operators on $L^2(C^*)$:
\begin{align*}
({\cal F}_{C^*} f)(\xi) &= -\frac1{\pi} \int_{\xi' \in C^*} \Psi_0\bigl(
\xi_1\xi'_1+\xi_2\xi'_2+\xi_3\xi'_3+\xi_4\xi'_4\bigr) \cdot f(\xi')
\,\frac{dS}{\|\xi'\|}, \\
(\FC f)(\xi) &= -\frac1{\pi} \int_{\xi' \in C^*}
\Psi_0\bigl(-\langle \xi, \xi' \rangle\bigr) \cdot f(\xi')
\,\frac{dS}{\|\xi'\|}, \hspace{1in} f \in L^2(C^*),
\end{align*}
where
\begin{equation}  \label{Psi_0}
\Psi_0(t) = \begin{cases}
Y_0(2\sqrt{2t}), & t>0;\\
-\frac2{\pi} K_0(2\sqrt{-2t}), & t<0;
\end{cases}
\end{equation}
$Y_0$ is the Bessel function of the second kind,
$K_0$ is the modified Bessel function of the second kind.
As a special case of Theorem 5.1.1 in \cite{KobM}, in the Schr\"odinger model
$(\pi, L^2(C^*))$ for the minimal representation of $O(3,3)$,
the operator $\pi(w_0): L^2(C^*) \to L^2(C^*)$ is given by
$\pi(w_0)f = {\cal F}_{C^*}f$.
Since $w'_0$ can be expressed as $w_0$ times a diagonal matrix with diagonal
entries $(-1,1,1,-1,-1,-1)$ and $-Id \in O(3,3)$ acts on the minimal
representation trivially, we immediately obtain that the operator
$\pi(w'_0): L^2(C^*) \to L^2(C^*)$ is given by $\pi(w'_0)f = \FC f$.

We finish this section by identifying
$\mathfrak{sp}(2,\BB R) \simeq \mathfrak{so}(3,2)$ sitting inside
$\mathfrak{sl}(2,\HR) \simeq \mathfrak{sl}(4,\BB R)$.
Set $\epsilon_1=\epsilon_2=1$, $\epsilon_3=\epsilon_4=-1$ and
let $\deg$ denote the degree operator plus identity:
$$
\deg f = \xi_1 \frac{\partial f}{\partial\xi_1}
+ \xi_2 \frac{\partial f}{\partial\xi_2}
+ \xi_3 \frac{\partial f}{\partial\xi_3}
+ \xi_4 \frac{\partial f}{\partial\xi_4} + f.
$$
The differential operators acting on $L^2(C^*)$ as
\begin{align*}
X_{jk} = &\epsilon_j \epsilon_k \xi_j \frac {\partial}{\partial \xi_k} -
\xi_k \frac {\partial}{\partial \xi_j}, \qquad 1 \le j < k \le 4, \\
&i(\epsilon_jP_j + 4\xi_j), \qquad 1 \le j \le 4, \qquad \text{where }
P_j= \epsilon_j \xi_j \square_{2,2} - 2\deg\circ\frac{\partial}{\partial\xi_j},
\end{align*}
generate a Lie algebra isomorphic to
$\mathfrak{sp}(2,\BB R) \simeq \mathfrak{so}(3,2)$.
We use Lemma 17 from \cite{FL} describing $\mathfrak{sl}(2,\HR)$ action
on the space of harmonic functions (see also its corrected version
at the beginning of Subsection 3.2 in \cite{FL2}).
By direct computation we can see that the
elements $\begin{pmatrix} A & 0 \\ 0 & D \end{pmatrix}$
with $A, D \in \HR$, $A^+=-A$, $D^+=-D$ act by linear combinations of
$X_{jk}$ which generate a Lie algebra isomorphic to
$\mathfrak{sl}(2,\BB R) \times \mathfrak{sl}(2,\BB R)$.
Next we see that elements $\begin{pmatrix} 0 & B \\ 0 & 0 \end{pmatrix}$
correspond to linear combinations of multiplication operators
$i\xi_j$, $1 \le j \le 4$. More specifically, choosing $B$ equal
$e_0$, $\tilde e_1$, $\tilde e_2$ and $e_3$ gives rise to multiplication by
$i\xi_1$, $i\xi_3$, $i\xi_4$ and $i\xi_2$ respectively.
Finally, by direct computation we get
\begin{multline*}
2 ( X \partial X - X) =
e_0 \Bigl( -N(X) \frac {\partial}{\partial x_1} +2x_1 \deg \Bigr)
+\tilde e_1 \Bigl( N(X) \frac{\partial}{\partial x_3} +2x_3\deg \Bigr) \\
+\tilde e_2 \Bigl( N(X) \frac{\partial}{\partial x_4} +2x_4\deg \Bigr)
+ e_3 \Bigl( -N(X) \frac {\partial}{\partial x_2} +2x_2 \deg \Bigr),
\end{multline*}
which implies that elements $\begin{pmatrix} 0 & 0 \\ C & 0 \end{pmatrix}$
correspond to linear combinations of operators
$$
iP_j(-1) =
i\epsilon_j \xi_j \square_{2,2} - i(2\deg+4)\circ\frac{\partial}{\partial\xi_j},
\qquad 1 \le j \le 4.
$$
More specifically, choosing $C$ equal $e_0$, $\tilde e_1$, $\tilde e_2$
and $e_3$ gives rise to operators $-iP_1(-1)$, $-iP_3(-1)$, $-iP_4(-1)$
and $iP_2(-1)$ respectively.
This implies that the Lie algebra generated by the differential operators
$X_{jk}$ and $i(\epsilon_jP_j + 4\xi_j)$ corresponds to
$$
\biggl\{ \begin{pmatrix} A & 4B \\ -B^+ & D \end{pmatrix} ;\: A,B,D \in \HR ,\:
A^+=-A,\: D^+=-D \biggr\} \subset \mathfrak{sl}(2,\HR).
$$
The ``4'' appears in the description of the Lie algebra because of the ``4''
in $i(\epsilon_jP_j + 4\xi_j)$.
The connected Lie subgroup of $GL(2,\HR)$ with the above Lie algebra
preserves the hyperboloid $\{ X \in \HR ;\: N(X)=-4 \}$
(see Lemma 25 in \cite{FL2}).

\section{$K$-finite Vectors of the Minimal Representation of $O(3,3)$}

In this section we find the $K$-finite vectors of the minimal representation
of $O(3,3)$ realized in $L^2(C^*)$.
This section is not needed for the proof of Theorem \ref{cont_ser_proj-2}.
Following \cite{KobM}, we renormalize the K-Bessel functions as
$$
\tilde K_n (r) = 2^n r^{-n} K_n (r), \qquad n \in \BB Z.
$$
We have the following relation between K-Bessel functions:
$$
K_{n+1}(r) - K_{n-1}(r) = 2n r^{-1} K_n(r),
$$
which implies
\begin{equation}  \label{K-rel}
r^2 \tilde K_{n+1}(2r) = n \tilde K_n(2r) + \tilde K_{n-1}(2r).
\end{equation}
The derivatives of K-Bessel functions satisfy the following relation:
$$
\Bigl( -\frac {d}{rdr} \Bigr)^m \bigl( r^{-n} K_n (r) \bigr) = r^{-n-m} K_{n+m} (r),
\qquad n \in \BB Z,
$$
which implies
$$
\Bigl( - \frac {2d}{rdr} \Bigr)^m \tilde K_n (r) = \tilde K_{n+m} (r),
\qquad
\frac d{dr} \tilde K_n (2r) = -2r \tilde K_{n+1} (2r),
\qquad n \in \BB Z.
$$

We identify $C^* \setminus \{0\}$ with $\BB R_+ \times S^1 \times S^1$
using coordinates
$$
(r, \theta_1, \theta_2) \:\longleftrightarrow\:
\xi=r(\cos\theta_1,\sin\theta_1,\cos\theta_2,\sin\theta_2) \quad \in C^*,
$$
which induces an isomorphism of Hilbert spaces
\begin{equation} \label{Hilbert-space-iso}
L^2(C^*) \simeq L^2 \Bigl( \BB R_+, \frac r2 dr \Bigr) \hat\otimes
L^2(S^1) \hat\otimes L^2(S^1).
\end{equation}
Consider the following functions on the cone $C^*$:
\begin{equation}  \label{K-vectors}
\tilde K_n(2r) (\xi_1 \pm i \xi_2)^l (\xi_3 \pm i \xi_4)^k
=
r^{l+k} \tilde K_n(2r) e^{\pm il\theta_1} e^{\pm ik\theta_2},
\end{equation}
$n \in \BB Z$, $k,l = 0,1,2,\dots$.
Clearly, they are $SO(2) \times SO(2)$-finite.

The action of $\mathfrak{sl}(4, \BB R)$ is generated by
multiplication by $\xi_i$, $1 \le i \le 4$, and differential operators
$P_j= \epsilon_j \xi_j \square_{2,2} - 2\deg\circ\frac{\partial}{\partial\xi_j}$,
$1 \le j \le 4$.
Hence it is important to find actions of these operators.
The multiplication operators are easy. For example, using (\ref{K-rel}):
\begin{multline*}
2\xi_1 \tilde K_n(2r) (\xi_1 \pm i \xi_2)^l
=
\tilde K_n(2r) (\xi_1 \pm i \xi_2)^{l+1}
+ r^2 \tilde K_n (2r) (\xi_1 \pm i \xi_2)^{l-1}  \\
=
\tilde K_n(2r) (\xi_1 \pm i \xi_2)^{l+1}
+ \bigl( (n-1) \tilde K_{n-1} (2r) + \tilde K_{n-2}(2r) \bigr)
(\xi_1 \pm i \xi_2)^{l-1};
\end{multline*}
\begin{multline*}
2i\xi_2 \tilde K_n(2r) (\xi_1 \pm i \xi_2)^l
=
\pm \tilde K_n(2r) (\xi_1 \pm i \xi_2)^{l+1}
\mp r^2 \tilde K_n (2r) (\xi_1 \pm i \xi_2)^{l-1}  \\
=
\pm \tilde K_n(2r) (\xi_1 \pm i \xi_2)^{l+1}
\mp \bigl( (n-1) \tilde K_{n-1} (2r) + \tilde K_{n-2}(2r) \bigr)
(\xi_1 \pm i \xi_2)^{l-1}.
\end{multline*}

To find the action of $P_j$'s we observe that in bipolar coordinates
$$
\square_{2,2} =
\frac{\partial^2}{\partial r_1^2} + \frac 1{r_1} \frac{\partial}{\partial r_1}
+ \frac 1{r_1^2} \frac{\partial^2}{\partial \theta_1^2}
-
\frac{\partial^2}{\partial r_2^2} - \frac 1{r_2} \frac{\partial}{\partial r_2}
- \frac 1{r_2^2} \frac{\partial^2}{\partial \theta_2^2}
$$
and
\begin{multline*}
- \Bigl(
\frac{\partial^2}{\partial r_2^2} + \frac 1{r_2} \frac{\partial}{\partial r_2}
+ \frac 1{r_2^2} \frac{\partial^2}{\partial \theta_2^2} \Bigr)
\bigl( \tilde K_n(2r_2) (\xi_3 \pm i \xi_4)^k \bigr) \\
=
4 \bigl( (k+1) \tilde K_{n+1}(2r_2)
- r_2^2 \tilde K_{n+2}(2r_2) \bigr) (\xi_3 \pm i \xi_4)^k \\
=
4\bigl((k-n) \tilde K_{n+1}(2r_2) - \tilde K_n(2r_2)\bigr) (\xi_3 \pm i\xi_4)^k.
\end{multline*}
We can realize  $r^{l+k} \tilde K_n(2r) e^{\pm il\theta_1} e^{\pm ik\theta_2}$
as the restriction of
$\tilde K_n(2r_2)(\xi_1 \pm i \xi_2)^l(\xi_3 \pm i \xi_4)^k$ to the cone $C^*$.
Then
\begin{multline*}
P_1 \bigl( \tilde K_n(2r_2) (\xi_1 \pm i \xi_2)^l (\xi_3 \pm i \xi_4)^k \bigr) \\
= - \biggl( \xi_1 \Bigl(
\frac{\partial^2}{\partial r_2^2} + \frac 1{r_2} \frac{\partial}{\partial r_2}
+ \frac 1{r_2^2} \frac{\partial^2}{\partial \theta_2^2} \Bigr) +2l\deg \biggr)
\bigl(\tilde K_n(2r_2) (\xi_1 \pm i\xi_2)^{l-1} (\xi_3 \pm i\xi_4)^k\bigr)\\
=
4 \xi_1 \bigl( (k-n) \tilde K_{n+1}(2r_2) - \tilde K_n(2r_2) \bigr)
(\xi_1 \pm i \xi_2)^l (\xi_3 \pm i \xi_4)^k \\
- 2l \bigl( (l+k) \tilde K_n(2r_2) - 2r_2^2 \tilde K_{n+1}(2r_2)\bigr)
(\xi_1 \pm i \xi_2)^{l-1} (\xi_3 \pm i \xi_4)^k  \\
=
2\bigl( (k-n) \tilde K_{n+1}(2r_2) - \tilde K_n(2r_2) \bigr)
(\xi_1 \pm i \xi_2)^{l+1} (\xi_3 \pm i \xi_4)^k \\
+ 2\bigl( n(k-n) \tilde K_n(2r_2) + (k-2n+1) \tilde K_{n-1}(2r_2)
- \tilde K_{n-2}(2r_2) \bigr) (\xi_1 \pm i \xi_2)^{l-1} (\xi_3 \pm i \xi_4)^k \\
+2l \bigl( (2n-l-k) \tilde K_n(2r_2)
+ 2 \tilde K_{n-1}(2r_2) \bigr) (\xi_1 \pm i \xi_2)^{l-1} (\xi_3 \pm i \xi_4)^k  \\
=
2 \bigl( (k-n)\tilde K_{n+1}(2r_2)
- \tilde K_n(2r_2) \bigr) (\xi_1 \pm i \xi_2)^{l+1} (\xi_3 \pm i \xi_4)^k \\
+ 2 \bigl( (n-l)(l+k-n)\tilde K_n(2r_2) + (2l+k-2n+1)\tilde K_{n-1}(2r_2)
- \tilde K_{n-2}(2r_2) \bigr) (\xi_1 \pm i \xi_2)^{l-1} (\xi_3 \pm i \xi_4)^k;
\end{multline*}
\begin{multline*}
P_2 \bigl( \tilde K_n(2r_2) (\xi_1 \pm i \xi_2)^l (\xi_3 \pm i \xi_4)^k \bigr)  \\
= -\biggl( \xi_2 \Bigl( \frac{\partial^2}{\partial r_2^2}
+ \frac 1{r_2} \frac{\partial}{\partial r_2}
+ \frac1{r_2^2} \frac{\partial^2}{\partial \theta_2^2} \Bigr)
\pm 2li \deg \biggr)
\bigl( \tilde K_n(2r_2) (\xi_1 \pm i \xi_2)^{l-1} (\xi_3 \pm i \xi_4)^k \bigr) \\
=
4 \xi_2 \bigl( (k-n) \tilde K_{n+1}(2r_2) - \tilde K_n(2r_2) \bigr)
(\xi_1 \pm i \xi_2)^l (\xi_3 \pm i \xi_4)^k \\
\mp 2li \bigl( (l+k) \tilde K_n(2r_2) - 2 r_2^2 \tilde K_{n+1}(2r_2) \bigr)
(\xi_1 \pm i \xi_2)^{l-1} (\xi_3 \pm i \xi_4)^k  \\
=
\mp 2i \bigl( (k-n) \tilde K_{n+1}(2r_2) - \tilde K_n(2r_2) \bigr)
(\xi_1 \pm i \xi_2)^{l+1} (\xi_3 \pm i \xi_4)^k \\
\pm 2i \bigl( n(k-n) \tilde K_n(2r_2) + (k-2n+1) \tilde K_{n-1}(2r_2)
- \tilde K_{n-2}(2r_2) \bigr) (\xi_1 \pm i \xi_2)^{l-1} (\xi_3 \pm i \xi_4)^k \\
\pm 2li \bigl( (2n-l-k) \tilde K_n(2r_2)
+ 2 \tilde K_{n-1}(2r_2) \bigr) (\xi_1 \pm i \xi_2)^{l-1} (\xi_3 \pm i \xi_4)^k  \\
=
\mp 2i \bigl( (k-n)\tilde K_{n+1}(2r_2)
- \tilde K_n(2r_2) \bigr) (\xi_1 \pm i \xi_2)^{l+1} (\xi_3 \pm i \xi_4)^k \\
\pm 2i \bigl( (n-l)(l+k-n)\tilde K_n(2r_2) + (2l+k-2n+1)\tilde K_{n-1}(2r_2)
- \tilde K_{n-2}(2r_2) \bigr) (\xi_1 \pm i \xi_2)^{l-1} (\xi_3 \pm i \xi_4)^k.
\end{multline*}
These calculations prove that $\mathfrak{sl}(4,\BB R)$ preserves the space
of finite linear combinations of the functions (\ref{K-vectors}).

\begin{prop}
The functions $\tilde K_n(2r) (\xi_1 \pm i \xi_2)^l (\xi_3 \pm i \xi_4)^k$,
$k, l \ge 0$ and $-\infty < n \le \min(k,l)$, are the
$K$-finite vectors for the action of $\mathfrak{sl}(4, \BB R)$.
Here we identify $\mathfrak{sl}(4, \BB R)$ with $\mathfrak{so}(3,3)$
and take the maximal compact subgroup corresponding to
$\mathfrak{so}(3) \times \mathfrak{so}(3)$.
\end{prop}

\pf
The Lie algebra $\mathfrak{so}(3) \times \mathfrak{so}(3)$ is generated
by $\mathfrak{so}(2) \times \mathfrak{so}(2)$ and elements
$N_j + \epsilon_j \B{N_j}$, $1 \le j \le 4$.
Then $\mathfrak{so}(2) \times \mathfrak{so}(2)$ acts by linear combinations
of $\partial / \partial \theta_1$ and $\partial / \partial \theta_2$,
and $N_j + \epsilon_j \B{N_j}$ acts by $2i\xi_j+ \frac i2 P_j$, $1 \le j \le 4$.
\begin{multline*}
\Bigl( 2(\xi_1+i\xi_2)+ \frac 12 (P_1 + i P_2) \Bigr)
\bigl( \tilde K_n(2r) (\xi_1 + i \xi_2)^l (\xi_3 \pm i \xi_4)^k \bigr)  \\
= 2 (k-n) \tilde K_{n+1}(2r) (\xi_1 + i \xi_2)^{l+1} (\xi_3 \pm i \xi_4)^k,
\end{multline*}
\begin{multline*}
\Bigl( 2(\xi_1+i\xi_2)+ \frac 12 (P_1 + i P_2) \Bigr)
\bigl( \tilde K_n(2r) (\xi_1 - i \xi_2)^l (\xi_3 \pm i \xi_4)^k \bigr)  \\
= 2 \bigl( (n-l)(l+k-n) \tilde K_n(2r) + (2l+k-n) \tilde K_{n-1}(2r) \bigr)
(\xi_1 - i \xi_2)^{l-1} (\xi_3 \pm i \xi_4)^k,
\end{multline*}
\begin{multline*}
\Bigl( 2(\xi_1-i\xi_2)+ \frac 12 (P_1 - i P_2) \Bigr)
\bigl( \tilde K_n(2r) (\xi_1 + i \xi_2)^l (\xi_3 \pm i \xi_4)^k \bigr)  \\
= 2 \bigl( (n-l)(l+k-n) \tilde K_n(2r) + (2l+k-n) \tilde K_{n-1}(2r) \bigr)
(\xi_1 + i \xi_2)^{l-1} (\xi_3 \pm i \xi_4)^k,
\end{multline*}
\begin{multline*}
\Bigl( 2(\xi_1-i\xi_2)+ \frac 12 (P_1 - i P_2) \Bigr)
\bigl( \tilde K_n(2r) (\xi_1 - i \xi_2)^l (\xi_3 \pm i \xi_4)^k \bigr)  \\
= 2 (k-n) \tilde K_{n+1}(2r) (\xi_1 - i \xi_2)^{l+1} (\xi_3 \pm i \xi_4)^k.
\end{multline*}
Thus $\mathfrak{so}(3) \times \mathfrak{so}(3)$ acts on a function of
type (\ref{K-vectors}) finitely if $n \le \min(k,l)$.
Lemma 3.4.1 in \cite{KobM} implies that these functions belong to $L^2(C^*)$.

We get all the $K$-finite vectors this way because setting $n=\min(k,l)$
we obtain the highest weight vector from each $K$-type of $L^2(C^*)$ by
Theorem 3.1.1 in \cite{KobM}.
\qed

\section{Fourier Transforms of $\frac 1{(N(X) \pm R^2 \pm i0)^2}$}

As a preparation for the proof of Theorem \ref{cont_ser_proj-2},
in this section we compute the Fourier transforms of the distributions
$\bigl( N(X) \pm R^2 \pm i0 \bigr)^{-2}$ on $\HR$.
The following proposition is closely related to the formulas from \S 2.8,
Chapter III of \cite{GS}, which describe the Fourier transforms of
the generalized functions $(R^2+Q \pm i0)^{\lambda}$ associated to a quadratic
form $Q$ on $\BB R^n$, where $\re \lambda < -n/2$.
However, in our case $n=4$ and $\lambda=-2=-n/2$, and we can think of it
as a ``limit'' of the setting of \cite{GS}.
Recall that the Fourier transform is normalized as in (\ref{FT-normalization}),
and let $J_0$ denote the Bessel function of the first kind.

\begin{prop}  \label{ft-computation}
Let $R>0$, then
$$
\lim_{\epsilon \to 0^+} \frac1{4\pi^2} \int_{X \in \HR}
\frac {e^{i\xi \cdot X} \,dV}{\bigl( N(X)-R^2 \pm i\epsilon \bigr)^2}
= \begin{cases}
\frac{\pi}4 Y_0 \bigl( R\sqrt{\langle \xi,\xi \rangle} \bigr)
\pm i\frac{\pi}4 J_0 \bigl( R\sqrt{\langle \xi,\xi \rangle} \bigr)
& \langle \xi,\xi \rangle>0;\\
-\frac12 K_0 \bigl( R\sqrt{-\langle \xi,\xi \rangle} \bigr)
& \langle \xi,\xi \rangle<0;
\end{cases}
$$
$$
\lim_{\epsilon \to 0^+} \frac1{4\pi^2} \int_{X \in \HR}
\frac {e^{i\xi \cdot X} \,dV}{\bigl( N(X)+R^2 \pm i\epsilon \bigr)^2}
= \begin{cases}
-\frac12 K_0 \bigl( R\sqrt{\langle \xi,\xi \rangle} \bigr)
& \langle \xi,\xi \rangle>0; \\
\frac{\pi}4 Y_0 \bigl( R\sqrt{-\langle \xi,\xi \rangle} \bigr)
\mp i\frac{\pi}4 J_0 \bigl( R\sqrt{-\langle \xi,\xi \rangle} \bigr)
& \langle \xi,\xi \rangle<0.
\end{cases}
$$
\end{prop}

\pf
We use the following integral representations of the Bessel functions
that can be found, for instance, in \cite{Ba} or \cite {GR}:
\begin{equation} \label{JY}
J_0(u) = \frac 2{\pi} \int_0^{\infty} \sin (u \cosh t) \,dt, \qquad
Y_0(u) = -\frac 2{\pi} \int_0^{\infty} \cos (u \cosh t) \,dt,
\end{equation}
\begin{equation} \label{K}
K_0(u) = \int_0^{\infty} \cos (u \sinh t) \,dt
= \int_0^{\infty} \exp (-u \cosh t) \,dt.
\end{equation}

Write $\xi=(\xi_1,\xi_2,\xi_3,\xi_4)$ and let
$r_1=\sqrt{ (\xi_1)^2 + (\xi_2)^2}$, $r_2=\sqrt{ (\xi_3)^2 + (\xi_4)^2}$.
Rotating $\xi$ if necessary, without loss of generality we can assume that
$\xi=(r_1,0,r_2,0)$. Using the formulas
$$
\int \frac{dx}{(a^2-x^2)^2} = \frac x{2a^2(a^2-x^2)}
+ \frac 1{4a^3} \log\Bigl(\frac{x+a}{x-a}\Bigr) +C, \qquad
\int \frac {dx}{(a^2+x^2)^{3/2}} = \frac x{a^2 \sqrt{a^2+x^2}} +C,
$$
we compute:
\begin{multline}  \label{ft-eqn}
\lim_{\epsilon \to 0^+} \iiiint_{\BB R^4} \frac {e^{i(r_1x_1+r_2x_3)}}
{( x_1^2+x_2^2-x_3^2-x_4^2-R^2+i\epsilon)^2} \,dx_1dx_2dx_3dx_4  \\
= -\frac{\pi i}2 \lim_{\epsilon \to 0^+} \iiint_{\BB R^3}
\frac {e^{i(r_1x_1+r_2x_3)}}{(x_1^2+x_2^2-x_3^2-R^2+i\epsilon)^{3/2}}\,dx_1dx_2dx_3 \\
= -\pi i \lim_{\epsilon \to 0^+} \iint_{\BB R^2}
\frac {e^{i(r_1x_1+r_2x_3)}}{x_1^2-x_3^2-R^2+i\epsilon} \,dx_1dx_3.
\end{multline}
Then we divide $\BB R^2$ into two regions: $\{|x_1|>|x_3|\}$ and 
$\{|x_1|<|x_3|\}$. On the first region we change the coordinates
$(x_1,x_3)$ to $(s,t)$ so that
$$
x_1 = Rs \cosh t, \quad x_3 = Rs \sinh t, \qquad -\infty < s < \infty,
\: -\infty < t < \infty,
$$
and on the second region we change the coordinates $(x_1,x_3)$ to $(s,t)$
so that
$$
x_1 = Rs \sinh t, \quad x_3 = Rs \cosh t, \qquad -\infty < s < \infty,
\: -\infty < t < \infty,
$$
and (\ref{ft-eqn}) reduces to
\begin{equation} \label{2-parts}
\pi i \lim_{\epsilon \to 0^+} \iint_{\BB R^2} |s| \biggl(
\frac{\exp\bigl(iRs(r_1\cosh t+r_2\sinh t)\bigr)}{1-s^2-i\epsilon} +
\frac{\exp\bigl(iRs(r_1\sinh t+r_2\cosh t)\bigr)}{1+s^2-i\epsilon}\biggr)\,dsdt.
\end{equation}

If $r_1>r_2$, then, for some $\tau \ge 0$, we can write
$$
r_1 = \sqrt{r_1^2-r_2^2} \cdot \cosh \tau, \qquad
r_2 = \sqrt{r_1^2-r_2^2} \cdot \sinh \tau,
$$
and
\begin{multline*}
\int_{\BB R} \exp \bigl(iRs(r_1\sinh t+r_2\cosh t)\bigr)\,dt
=
\int_{\BB R} \exp \Bigl( iRs \sqrt{r_1^2-r_2^2} \sinh(t+\tau) \Bigr)\,dt  \\
=
\int_{\BB R} \cos \Bigl( Rs \sqrt{r_1^2-r_2^2} \sinh(t+\tau) \Bigr)\,dt
+ i \int_{\BB R} \sin \Bigl( Rs \sqrt{r_1^2-r_2^2} \sinh(t+\tau) \Bigr)\,dt\\
= 2 K_0 \Bigl( Rs \sqrt{r_1^2-r_2^2} \Bigr)
= \int_{\BB R} \exp \Bigl( -R|s| \sqrt{r_1^2-r_2^2} \cosh(t+\tau) \Bigr)\,dt \\
= \int_{\BB R} \exp \bigl(-R|s|(r_1\cosh t+r_2\sinh t)\bigr)\,dt.
\end{multline*}
Thus we can rewrite (\ref{2-parts}) as
\begin{multline*}
\pi i \lim_{\epsilon \to 0^+} \int_{\BB R} \biggl( \int_{\Gamma_1}
\frac{s \exp \bigl(iRs(r_1\cosh t+r_2\sinh t) \bigr)}{1-s^2-i\epsilon}\,ds \\
+ \int_{\Gamma_2} \frac{s \exp \bigl(-iRs(r_1\cosh t+r_2\sinh t) \bigr)}
{1-s^2-i\epsilon}\,ds \biggr)\,dt,
\end{multline*}
where $\Gamma_1$ is the contour in the complex plane $\BB C$ which starts at
infinity, goes along the positive imaginary axis towards the origin and
then along the positive real axis towards infinity, and $\Gamma_2$ is the
contour which starts at infinity, goes along the negative imaginary axis
towards the origin and then along the positive real axis towards infinity.
We can compute the integrals with respect to $s$ using residues:
\begin{multline*}
\lim_{\epsilon \to 0^+} \int_{\HR}
\frac{e^{i(r_1x_1+r_2x_3)}}{\bigl( N(X)-R^2+i\epsilon \bigr)^2} \,dV
= -\pi^2 \int_{\BB R} \exp \bigl(-iR(r_1\cosh t+r_2\sinh t) \bigr)\,dt\\
= -\pi^2 \int_{\BB R} \cos \Bigl( R \sqrt{r_1^2-r_2^2} \cosh(t+\tau) \Bigr)\,dt
+ i\pi^2 \int_{\BB R} \sin \Bigl( R \sqrt{r_1^2-r_2^2} \cosh(t+\tau) \Bigr)\,dt \\
= \pi^3 Y_0 \Bigl( R\sqrt{r_1^2-r_2^2} \Bigr)
+i\pi^3 J_0 \Bigl( R\sqrt{r_1^2-r_2^2} \Bigr).
\end{multline*}

If $r_1<r_2$, then, for some $\tau \ge 0$, we can write
$$
r_1 = \sqrt{r_2^2-r_1^2} \cdot \sinh \tau, \qquad
r_2 = \sqrt{r_2^2-r_1^2} \cdot \cosh \tau,
$$
and
\begin{multline*}
\int_{\BB R} \exp \bigl(iRs(r_1\cosh t+r_2\sinh t) \bigr)\,dt
=
\int_{\BB R} \exp \Bigl( iRs \sqrt{r_2^2-r_1^2} \sinh(t+\tau) \Bigr)\,dt  \\
=
\int_{\BB R} \cos \Bigl( Rs \sqrt{r_2^2-r_1^2} \sinh(t+\tau) \Bigr)\,dt
+ i \int_{\BB R} \sin \Bigl( Rs \sqrt{r_2^2-r_1^2} \sinh(t+\tau) \Bigr)\,dt\\
= 2 K_0 \Bigl( Rs \sqrt{r_2^2-r_1^2} \Bigr)
= \int_{\BB R} \exp \Bigl( -R|s| \sqrt{r_2^2-r_1^2} \cosh(t+\tau) \Bigr)\,dt \\
= \int_{\BB R} \exp \bigl(-R|s|(r_1\sinh t+r_2\cosh t)\bigr)\,dt.
\end{multline*}
Thus we can rewrite (\ref{2-parts}) as
\begin{multline*}
\pi i \lim_{\epsilon \to 0^+} \int_{\BB R} \biggl( \int_{\Gamma_1}
\frac {s \exp \bigl( iRs(r_1\sinh t+r_2\cosh t) \bigr)}{1+s^2-i\epsilon}\,ds \\
+ \int_{\Gamma_2} \frac {s \exp \bigl(-iRs(r_1\sinh t+r_2\cosh t) \bigr)}
{1+s^2-i\epsilon}\,ds \biggr)\,dt,
\end{multline*}
where $\Gamma_1$ and $\Gamma_2$ are the same contours in the complex plane
$\BB C$ as above.
Using residues to compute the integrals with respect to $s$ we get:
\begin{multline*}
\lim_{\epsilon \to 0^+} \int_{\HR}
\frac{e^{i(r_1x_1+r_2x_3)}}{\bigl( N(X)-R^2+i\epsilon \bigr)^2} \,dV
= -\pi^2 \int_{\BB R} \exp \bigl(-R(r_1\sinh t+r_2\cosh t) \bigr)\,dt \\
= -\pi^2 \int_{\BB R} \exp \bigl(-R \sqrt{r_2^2-r_1^2} \cosh(t+\tau)) \bigr)\,dt
= -2\pi^2 K_0 \Bigl( R\sqrt{r_2^2-r_1^2} \Bigr).
\end{multline*}
The other Fourier transforms are found in exactly the same way.
\qed

\begin{cor}
Let $R>0$ and $\xi, \: \xi' \in C^*$, then
$$
\lim_{\epsilon \to 0^+} \frac1{4\pi^2} \int_{X \in \HR} \biggl(
\frac{e^{i(\xi-\xi') \cdot X}}{\bigl( N(X)-4 + i\epsilon \bigr)^2}
+ \frac{e^{i(\xi-\xi') \cdot X}}{\bigl( N(X)-4 - i\epsilon \bigr)^2} \biggr)\,dV
= \frac{\pi}2 \cdot \Psi_0 \bigl( -\langle \xi, \xi' \rangle \bigr),
$$
\begin{multline*}
\lim_{\epsilon \to 0^+} \frac1{4\pi^2} \int_{X \in \HR} \biggl(
\frac{e^{i(\xi-\xi') \cdot X}}{\bigl( N(X)-R^2 + i\epsilon \bigr)^2}
- \frac{e^{i(\xi-\xi') \cdot X}}{\bigl( N(X)-R^2 - i\epsilon \bigr)^2} \biggr)\,dV \\
= \begin{cases}
0 & \langle \xi,\xi' \rangle>0, \\
\frac{\pi i}2 J_0 \bigl( R\sqrt{-2\langle \xi,\xi' \rangle} \bigr)
& \langle \xi,\xi' \rangle <0
\end{cases}
= \frac{\pi i}2 \cdot \Phi^+_0\Bigl(-\frac{R^2}4 \langle\xi,\xi'\rangle\Bigr).
\end{multline*}
\end{cor}

\pf
Since $\xi, \xi' \in C^*$,
$\langle \xi-\xi', \xi-\xi' \rangle = -2 \langle \xi, \xi' \rangle$,
and the result follows from the above proposition.
\qed

We conclude this section with the following expressions for the kernels
$\Psi_0 \bigl( \langle \xi, \xi' \rangle \bigr)$ and
$\Phi^+_0 \bigl( \langle \xi, \xi' \rangle \bigr)$
defined by (\ref{Psi_0}) and (\ref{Phi_0}).

\begin{lem}
Let $R>0$,
$\xi=(\xi_1,\xi_2,\xi_3,\xi_3), \: \xi'=(\xi'_1,\xi'_2,\xi'_3,\xi'_3) \in C^*$,
$$
r_1=\sqrt{ (\xi_1-\xi'_1)^2 + (\xi_2-\xi'_2)^2}, \qquad
r_2=\sqrt{ (\xi_3-\xi'_3)^2 + (\xi_4-\xi'_4)^2},
$$
then
\begin{align*}
\Psi_0 \Bigl( \frac{R^2}4 \langle \xi, \xi' \rangle \Bigr)
&= -\frac1{\pi} \int_{-\infty}^{\infty} \cos(Rr_1 \sinh t + Rr_2 \cosh t) \,dt, \\
\Psi_0 \Bigl( -\frac{R^2}4 \langle \xi, \xi' \rangle \Bigr) 
&= -\frac1{\pi} \int_{-\infty}^{\infty} \cos(Rr_1 \cosh t + Rr_2 \sinh t) \,dt, \\
\Phi^+_0 \Bigl( \frac{R^2}4 \langle \xi, \xi' \rangle \Bigr)
&= \frac1{\pi} \int_{-\infty}^{\infty} \sin(Rr_1 \sinh t + Rr_2 \cosh t) \,dt, \\
\Phi^+_0 \Bigl( -\frac{R^2}4 \langle \xi, \xi' \rangle \Bigr) 
&= \frac1{\pi} \int_{-\infty}^{\infty} \sin(Rr_1 \cosh t + Rr_2 \sinh t) \,dt.
\end{align*}
\end{lem}

Since we are not going to use this lemma, we omit its proof.

\section{Fourier Transform of $\Pl'_R$ and
the Proof of Theorem \ref{cont_ser_proj-2}}

In this section we compute the Fourier transform of $\Pl'_R$,
compare it with $\FC$ and use that comparison to finish
the proof of Theorem \ref{cont_ser_proj-2}.
As in \cite{KobM}, let $T: L^2(C^*) \to {\cal S}'(\HR^*)$ be the continuous
injective map sending a function $f \in L^2(C^*)$ into a distribution
$f \cdot \delta(C^*)$ belonging to the space of tempered distributions on
$\HR^*$. If a function $\phi$ on $\HR$ satisfies $\square_{2,2}\phi=0$,
its Fourier transform $\hat\phi$ is a distribution supported on the cone $C^*$.
As explained in \cite{KobM}, \cite{KobO}, if this function $\phi$ lies in the
solution model for minimal representation $V^{3,3}$ of $O(3,3)$,
then $\hat\phi$ lies in the image $T(L^2(C^*))$, and so it makes sense to talk
about $T^{-1} (\hat\phi) \in L^2(C^*)$.
Let $FT: V^{3,3} \to L^2(C^*)$ denote the map $\phi \mapsto T^{-1} (\hat\phi)$.
Then we have a commutative diagram
$$
\begin{CD}
V^{3,3} @>{FT}>> L^2(C^*)\\
@V{\varpi^{3,3}(w'_0)}VV @VV{\FC}V \\
V^{3,3} @>{FT}>> L^2(C^*)
\end{CD}
$$
As a representation of $SU(1,1) \times SU(1,1)$, which is a double cover of
$SO(2,2)$, $L^2(C^*)$ is a direct integral of the spaces of homogeneous
functions of degree $-1+i\rho$ with $\rho \in \BB R$.
Then $FT^{-1} \bigl(L^2(C^*)\bigr)$ is also a direct integral of the spaces of
homogeneous functions of degree $-1+i\rho$, $\rho \in \BB R$, satisfying
$\square_{2,2}\phi=0$.
By \cite{St} these functions restrict to the continuous spectrum functions
on $SU(1,1)$ realized as the unit hyperboloid in $\HR$. This proves that
the minimal representation $V^{3,3}$ is nothing but the closure of the space
generated by the wave packets of matrix coefficients $t^l_{n\,\underline{m}}$'s
of $SU(1,1)$ with $\re l = -1/2$.

Since the Fourier transform replaces products with convolutions and vice versa,
we obtain the following integral expression for the operator
$\varpi^{3,3}(w'_0)$ in the solution model:
\begin{multline*}
\phi(X) \mapsto (\varpi^{3,3}(w'_0) \phi)(W) \\
= -\frac2{\pi^2} \lim_{\epsilon \to 0^+} \int_{X \in C} \biggl(
\frac1{\bigl(N(W-X)-4 +i\epsilon\bigr)^2} +
\frac1{\bigl(N(W-X)-4 -i\epsilon\bigr)^2} \biggr) \phi(W-X) \,\frac{dS}{\|X\|}.
\end{multline*}
Here we used that the Fourier transform of $\delta(C)$ is $\delta(C^*)$, where
$$
\delta(C^*): \quad \psi \mapsto
\frac 12 \int_{\xi \in C^*} \psi(\xi) \,\frac{dS}{\|\xi\|}
$$
(see \cite{GS}, Chapter III, Section 2.6).

Next we define a map $\widehat{\Pl}'_R$ on $L^2(C^*)$ by
$$
(\widehat{\Pl}'_R f)(\xi) = \frac{i}{4\pi} \int_{\xi' \in C^*}
 \Phi^+_0 \Bigl( -\frac{R^2}4 \langle \xi, \xi' \rangle \Bigr) \cdot f(\xi')
\,\frac{dS}{\|\xi'\|}, \qquad  f \in L^2(C^*),
$$
where
\begin{equation}  \label{Phi_0}
\Phi^+_0(t) =
\begin{cases} J_0 \bigl( 2\sqrt{2t} \bigr) & t>0; \\ 0 & t<0. \end{cases}
\end{equation}

\begin{lem}
The map $\widehat{\Pl}'_R$ sends $L^2(C^*)$ into $L^2(C^*)$.
\end{lem}

\pf
Recall the Hilbert space isomorphism (\ref{Hilbert-space-iso}), hence
$$
L^2(C^*) \simeq \widehat{\bigoplus}_{l,k \in \BB Z}
L^2\Bigl( \BB R_+, \frac{r}2 dr \Bigr)
\otimes \BB C e^{i\theta_1} \otimes \BB C e^{i\theta_2}.
$$
Note that the map $\widehat{\Pl}'_R$ is $S^1 \times S^1$-equivariant.
Let $(\widehat{\Pl}'_R)_{k,l}$ denote the restriction of
$\widehat{\Pl}'_R$ to
$L^2(\BB R_+, \frac{r}2 dr) \otimes \BB C e^{i\theta_1} \otimes \BB C e^{i\theta_2}$.
Thus, $(\widehat{\Pl}'_R)_{k,l}$ is essentially a map from functions on $\BB R_+$
into functions on $\BB R_+$, and, clearly, $(\widehat{\Pl}'_R)_{k,l}$ is
well-defined on smooth compactly supported functions on $\BB R_+$.
The fact that $(\widehat{\Pl}'_R)_{k,l}$ extends to $L^2(\BB R_+, \frac{r}2 dr)$
and that its image is contained in $L^2(\BB R_+, \frac{r}2 dr)$ is proven in
Section 5.4 of \cite{KobM} because $\Phi^+_0(t)$ is the integral kernel of
the operator ${\cal F}_{C^*}$ for the Minkowski space $\BB M$
(i.e. the case $p=4$, $q=2$).
\qed

The same reasoning as above shows that the map $\Pl'_R$ restricted to $V^{3,3}$
can be fit into a diagram
$$
\begin{CD}
V^{3,3} @>{FT}>> L^2(C^*)\\
@V{\Pl'_R}VV @VV{\widehat{\Pl}'_R}V \\
V^{3,3} @>{FT}>> L^2(C^*)
\end{CD}
$$

Note that if $\phi$ is a function on $\HR$ which is homogeneous of degree
$-1+i\rho$, then $FT\phi$ is homogeneous of degree $-1-i\rho$, and both
$\widehat{\Pl}'_R (FT\phi)$, $\FC(FT\phi)$ are homogeneous of degree $-1+i\rho$.
The most obvious thing to do is to apply operators $\widehat{\Pl}'_R$ and
$\FC$ to functions $f$ on $C^*$ that are homogeneous of degree $1-i\rho$,
$\rho \in \BB R$, then find the ratio $ (\widehat{\Pl}'_R f) / (\FC f)$
and compare it with (\ref{w_0-action}).
Unfortunately, such functions are not in $L^2(C^*)$ and we would run into
convergence problems.

Fix any $\xi \in C^* \setminus \{0\}$ and $\epsilon \in \{0, 1\}$.
We identify $C^*$ with $\BB R_+ \times S^1 \times S^1$.
Let $\psi$ be a smooth function on $S^1 \times S^1$ such that
$\psi(-\theta_1, -\theta_2) = (-1)^{\epsilon} \psi(\theta_1, \theta_2)$.
We shall regard $\psi$ as a function on $C^*$ and assume that $\psi$ satisfies
$$
\supp \psi \cap \{\xi' \in C^*; \langle \xi, \xi' \rangle =1 \}
$$
is compact.
We define a function $f_{\xi,\epsilon}$ on $C^*$ by
$$
f_{\xi,\epsilon}(\xi') = \psi \cdot |\langle \xi, \xi' \rangle|^{-1/2} \cdot
\exp(-|\langle \xi, \xi' \rangle|).
$$
By design, $f_{\xi,\epsilon}$ belongs to $L^2(C^*)$.
Next we compute $\widehat{\Pl}'_R f_{\xi,\epsilon}(s\xi)$ and
$\FC f_{\xi,\epsilon}(s\xi)$, $s \in \BB R$.
For this purpose we use the following integrals, which can be found in
\cite{GR} (special cases of 2.667(7-8)):
$$
\int_0^{\infty} t^2 e^{-at} \sin(bt)\,dt = \frac{2b(3a^2-b^2)}{(a^2+b^2)^3}, \qquad
\int_0^{\infty} t^2 e^{-at} \cos(bt)\,dt = \frac{2a(a^2-3b^2)}{(a^2+b^2)^3}, \qquad
a>0,
$$
$$
\text{and} \qquad
\int_0^{\infty} s^{\mu-1} e^{-as}\,ds = a^{-\mu} \Gamma(\mu).
$$
Using (\ref{JY}) we find that $\widehat{\Pl}'_R f_{\xi,\epsilon}(s\xi)$ is
proportional to the integral with respect to $\theta$
(from $\theta=0$ to $\infty$) of
\begin{multline*}
(-1)^{\epsilon} \frac{i}{2\pi^2}
\int_0^{\infty}\sqrt{t} e^{-\sqrt{t}} \sin(R\sqrt{2st} \cosh \theta) \,dt
= (-1)^{\epsilon} \frac{i}{\pi^2}
\int_0^{\infty} t^2 e^{-t} \sin(R\sqrt{2s} \cosh \theta t) \,dt \\
= (-1)^{\epsilon} \frac{2\sqrt{2}i}{\pi^2}
\frac{3R\cosh\theta s^{1/2} -2R^3\cosh^3\theta s^{3/2}}{(1+2R^2\cosh^2\theta s)^3}.
\end{multline*}
Similarly, using (\ref{JY}) and (\ref{K}) we find that 
$\FC f_{\xi,\epsilon}(s\xi)$ is proportional
to the integral with respect to $\theta$ (from $\theta=0$ to $\infty$) of
\begin{multline*}
\frac{2}{\pi^2}
\int_0^{\infty}\sqrt{t} e^{-\sqrt{t}} \Bigl( \exp(-2\sqrt{2st} \cosh \theta)
+ (-1)^{\epsilon} \cos(2\sqrt{2st} \cosh\theta) \Bigr) \,dt \\
= \frac4{\pi^2} \int_0^{\infty} t^2 e^{-t} \Bigl( \exp(-2\sqrt{2s} \cosh \theta t)
+ (-1)^{\epsilon} \cos(2\sqrt{2s} \cosh\theta t) \Bigr) \,dt \\
= \frac8{\pi^2} \biggl( \frac1{(1+2\sqrt{2s} \cosh\theta)^3}
+ (-1)^{\epsilon} \frac{1-24\cosh^2\theta s}{(1+8\cosh^2\theta s)^3} \biggr).
\end{multline*}

Next we apply the Mellin transform
$f(s) \mapsto ({\cal M}f)(\rho) = \int_0^{\infty} f(s) s^{1-i\rho}\,\frac{ds}s$.
For this purpose we use integral formula 8.384 from \cite{GR}
$$
\int_0^{\infty} \frac{s^{\mu-1}}{(1+as)^{\nu}}\,ds
= a^{-\mu} \frac{\Gamma(\mu)\Gamma(\nu-\mu)}{\Gamma(\nu)}.
$$
Thus ${\cal M} \bigl(\widehat{\Pl}'_R f_{\xi,\epsilon}(s\xi)\bigr)$ is
proportional to the integral with respect to $\theta$ of
\begin{multline*}
(-1)^{\epsilon} \frac{2\sqrt{2}i}{\pi^2} \int_0^{\infty}
\frac{3R\cosh\theta s^{1/2} -2R^3\cosh^3\theta s^{3/2}}
{(1+2R^2\cosh^2\theta s)^3} s^{-i\rho} \,ds\\
= (-1)^{\epsilon} \frac{i}{\pi^2}(\sqrt{2}R\cosh\theta)^{-2+2i\rho}
\bigl(3 \Gamma(3/2-i\rho)\Gamma(3/2+i\rho)
- \Gamma(5/2-i\rho)\Gamma(1/2+i\rho) \bigr)\\
= (-1)^{\epsilon} \frac{i}{\pi^2}(\sqrt{2}R\cosh\theta)^{-2+2i\rho}
\Gamma(1/2-i\rho)\Gamma(1/2+i\rho) (1/2-i\rho)
\bigl(3(1/2+i\rho) - (3/2-i\rho) \bigr) \\
= (-1)^{\epsilon+1} \frac2{\pi} \frac{\rho(1-2i\rho)}{\cos(\pi i \rho)}
(\sqrt{2}R\cosh\theta)^{-2+2i\rho}.
\end{multline*}
On the other hand, ${\cal M} \bigl(\FC f_{\xi,\epsilon}(s\xi)\bigr)$
is proportional to the integral with respect to $\theta$ of
\begin{multline*}
\frac8{\pi^2} \int_0^{\infty} \biggl( \frac1{(1+2\sqrt{2s} \cosh\theta)^3}
+ (-1)^{\epsilon} \frac{1-24\cosh^2\theta s}{(1+8\cosh^2\theta s)^3} \biggr)
s^{-i\rho} \,ds\\
= \frac4{\pi^2}(2\sqrt{2}\cosh\theta)^{-2+2i\rho}
\Bigl( 2\Gamma(2-2i\rho)\Gamma(1+2i\rho)
+ (-1)^{\epsilon} \bigl( \Gamma(1-i\rho)\Gamma(2+i\rho)
- 3\Gamma(2-i\rho)\Gamma(1+i\rho) \bigr) \Bigr)\\
= \frac{4i\rho}{\pi^2}(2\sqrt{2}\cosh\theta)^{-2+2i\rho}
\Bigl( 4(1-2i\rho)\Gamma(1-2i\rho)\Gamma(2i\rho)
+ (-1)^{\epsilon} \Gamma(1-i\rho)\Gamma(i\rho)
\bigl((1+i\rho) - 3(1-i\rho)\bigr) \Bigr)\\
= \frac{8i\rho(1-2i\rho)}{\pi}(2\sqrt{2}\cosh\theta)^{-2+2i\rho}
\Bigl( \frac2{\sin(2\pi i\rho)} - \frac{(-1)^{\epsilon}}{\sin(\pi i\rho)} \Bigr)\\
= \frac{8i}{\pi} \frac{\rho(1-2i\rho)}{\cos(\pi i\rho)}
(2\sqrt{2}\cosh\theta)^{-2+2i\rho} \cdot
\begin{cases} \tan(\pi i\rho/2) & \text{if $\epsilon=0$;}\\
\cot(\pi i\rho/2) & \text{if $\epsilon=1$.} \end{cases}
\end{multline*}
Therefore,
\begin{multline*}
\frac{{\cal M} \bigl(\widehat{\Pl}'_R f_{\xi,\epsilon}(s\xi)\bigr)}
{{\cal M} \bigl(\FC f_{\xi,\epsilon}(s\xi)\bigr)}(\rho)
= i R^{-2+2i\rho} 2^{-2i\rho} \cdot
\begin{cases} \cot(\pi i\rho/2) & \text{if $\epsilon=0$;}\\
-\tan(\pi i\rho/2) & \text{if $\epsilon=1$} \end{cases} \\
= R^{-2+2i\rho} 2^{-2i\rho} \cdot
\begin{cases} \coth(\pi\rho/2) & \text{if $\epsilon=0$;}\\
\tanh(\pi\rho/2) & \text{if $\epsilon=1$.} \end{cases}
\end{multline*}
This completes the proof of Theorem \ref{cont_ser_proj-2}.

\separate

\separate

\noindent
{\em Department of Mathematics, Yale University,
P.O. Box 208283, New Haven, CT 06520-8283}\\
{\em Department of Mathematics, Indiana University,
Rawles Hall, 831 East 3rd St, Bloomington, IN 47405}

\noindent
Corresponding author's e-mail: mlibine$@$indiana.edu

\end{document}